\documentstyle[amssymb,amsfonts]{amsart}

\def\R{{\hbox{\bf R}}}

 at 10 true pt

\def\allt#1{%
\smash{
\vtop{%
     \ialign{%
        ##\crcr
        $\hfil\displaystyle{\tilde \forall}\hfil$\crcr%
        \noalign{\kern1.5pt\nointerlineskip}
        $\hfil\!\!#1\hfil$\crcr\noalign{\kern1.5pt}
        }
       }
      } \hbox{$\vphantom{#1}$}
     }

\def\be#1{\begin{equation}\label{#1}}
\def\bas{\begin{align*}}
\def\eas{\end{align*}}
\def\bi{\begin{itemize}}
\def\ei{\end{itemize}}

\def\eps{\varepsilon}
\newenvironment{proof}{\noindent {\bf Proof} }{\endprf\par}
\def \endprf{\hfill  {\vrule height6pt width6pt depth0pt}\medskip}
\def\emph#1{{\it #1}}
\def\textbf#1{{\bf #1}}

\theoremstyle{plain}
  \newtheorem{theorem}[subsection]{Theorem}

  \newtheorem{lemma}[subsection]{Lemma}
  \newtheorem{corollary}[subsection]{Corollary}

\theoremstyle{remark}

\theoremstyle{definition}

\include{psfig}

\begin{document}

\title[Arithmetic projections and Kakeya]{Bounds on arithmetic projections, and applications to the Kakeya conjecture}

\author{Nets Hawk Katz}
\address{Department of Mathematics, University of Illinois at Chicago, Chicago IL 60607-7045}
\email{nets@@math.uic.edu}

\author{Terence Tao}
\address{Department of Mathematics, UCLA, Los Angeles CA 90095-1555}
\email{tao@@math.ucla.edu}


\begin{abstract}  Let $A$, $B$, be finite subsets of an
abelian group, and let $G \subset A \times B$ be such that
$\# A, \# B, \# \{ a+b: (a,b) \in G \} \leq N$.  We consider the 
question of estimating
the quantity $\# \{ a-b: (a,b) \in G \}$. 
In \cite{borg:high-dim} Bourgain obtained the bound
of $N^{2-\frac{1}{13}}$, and applied
this to the Kakeya conjecture.  We improve Bourgain's estimate
to $N^{2-\frac{1}{6}}$, and obtain the further improvement of $N^{2 - \frac{1}{4}}$ under the additional assumption $\# \{ a+2b: (a,b) \in G\} \leq N$.  As an application we conclude that Besicovitch sets
in $\R^n$ have Minkowski dimension at least
$\frac{4n}{7} + \frac{3}{7}$.  This is new for $n > 8$.
\end{abstract}

\maketitle

\section{Introduction}

Let $N$ be a positive integer, and let $(Z,+)$ be an abelian 
group.  Let $A$, $B$,
be finite subsets of $Z$ with cardinality
\be{ab-card}
\# A, \# B \leq N.
\end{equation}
Let $G$ be a subset of $A \times B$.  We consider the question of
bounding the quantity
\be{diff}
\# \{ a-b: (a,b) \in G \}.
\end{equation}
Without any further assumptions on $G$ one can only obtain the trivial bound of $N^2$.  However, in \cite{borg:high-dim} Bourgain showed that under the additional assumption
\be{c-card}
\# C \leq N, \hbox{ where } C = \{ a + b: (a,b) \in G \}
\end{equation}
that one could improve the bound on \eqref{diff} to $N^{2 - \frac{1}{13}}$.
The purpose of this paper is to obtain the following additional improvements.

\begin{theorem}\label{main}  
Let the notation and assumptions be as above.  Then we have
\be{diff-6}
\eqref{diff} \leq N^{2 - \frac{1}{6}}.
\end{equation}
If we make the further additional assumption
\be{d-card}
\# D \leq N, \hbox{ where } D = \{ a + 2b: (a,b) \in G \}
\end{equation}
then we may improve this further to 
\be{diff-4}
\eqref{diff} \leq N^{2 - \frac{1}{4}}.
\end{equation}
\end{theorem}

It appears plausible that further improvements are possible by adding further assumptions of the type in \eqref{c-card}, \eqref{d-card}, but our methods do not seem to be able take advantage of such additional assumptions.

By the arguments in \cite{borg:high-dim} this implies a bound on the 
Minkowski and Hausdorff dimensions of Besicovitch sets\footnote{To obtain the Minkowski bound, one adapts Proposition 1.7 of \cite{borg:high-dim} to the arguments in this paper with $E$ equal to the slices of the Besicovitch set at $x_n = 0,1,1/2,2/3$, and then averages over translations.}.  Recall that a 
Besicovitch set in $\R^n$, $n > 1$ is a set which contains a unit line segment
in every direction.  The Kakeya conjecture states that such sets must have
full dimension.

\begin{corollary}\cite{borg:high-dim} Let $E$ be a Besicovitch set in $\R^n$.  Then the Minkowski dimension of $E$ is at least $\frac{4n}{7} + 
\frac{3}{7}$, while the Hausdorff dimension of $E$ is at least $\frac{6n}{11} + \frac{5}{11}$.
\end{corollary}

The lower bound of $\frac{n+2}{2}$ for both types of dimension was obtained in
\cite{wolff:kakeya}.  Thus the Minkowski bound is new for $n>8$ and the Hausdorff bound is new for $n > 12$.  We remark that
Theorem \ref{main} can also be used to slightly improve some other
recent work on the Kakeya problem in \cite{borg:high-dim} and \cite{katzlabatao}.  Roughly speaking, the connection between Besicovitch sets and Theorem \ref{main} arises from taking $A$, $B$, $C$, $D$ to essentially be the slices of
the Besicovitch set at the hyperplanes $\{ x_n = 0\}$, $\{ x_n = 1\}$,
$\{ x_n = 1/2 \}$, and $\{ x_n = 2/3 \}$ respectively, and $G$ to be the set of
pairs in $A \times B$ whose associated line segment is contained in the Besicovitch set.  In principle the Hausdorff bound is improvable to match the Minkowski bound, but one must first prove an analogue of the results in
\cite{heath-brown}, namely that that subsets of $\{1, \ldots, N\}$, with density at least $1/(\log N)^\eps$ contain four distinct elements which are affinely equivalent to $\{ 0, 1, 1/2, 2/3 \}$ if $N$ is sufficiently large and $\eps$ sufficiently small.  For recent progress on this type of problem see Gowers \cite{gowers}.

Bourgain's bound of $N^{2 - \frac{1}{13}}$ was obtained using some ideas of Gowers \cite{gowers} in his work on the Balog-Szemer\'edi theorem
\cite{balog}.  Our arguments are more elementary, and will not yield any new results of Balog-Szemer\'edi type.

In the converse direction, a simple variant of an example in \cite{ruzsa}
shows that \eqref{diff} can be as large as 
$N^{\log(6)/\log(3)} = N^{2 - 0.36907\ldots}$ if one assumes \eqref{c-card},
or as large as $N^{\log(8)/\log(4)} = N^{2 - 0.5}$ if one also assumes \eqref{d-card}.  To see the former claim, let $n, M$ be large integers and
and for every integer $0 \leq a < M^n$
let $d_0(a), \ldots, d_{n-1}(a) \in \{0, \ldots, M-1 \}$ be the digits of $a$ 
base $M$, thus
$$ a = \sum_{i=0}^{n-1} d_i(a) M^i.$$
If we set $Z$ to be the integers and
\bas
A = B &= \{ 0 \leq a < M^n: d_i(a) \in \{ 0,1,3 \} \hbox{ for all }
i = 0, \ldots, n-1 \}\\
C &= \{ 0 \leq c < M^n: d_i(c) \in \{ 1,3,4 \} \hbox{ for all }
i = 0, \ldots, n-1 \}\\
G &= \{ (a,b) \in A \times B: d_i(a) \neq d_i(b) \hbox{ for all }
i = 0, \ldots, n-1 \}
\end{align*}
(cf. \cite{ruzsa}) then we easily see that the hypotheses are satisfied if $M \geq 7$ with $N = 3^n$ and $\eqref{diff} = \# G = 6^n$, hence the claim.  The latter claim is similar but uses the sets
\bas
A &= \{ 0 \leq a < M^n: d_i(a) \in \{ 0,2,3,4 \} \hbox{ for all }i = 0, \ldots, n-1 \}\\
B &= \{ 0 \leq b < M^n: d_i(b) \in \{ 0,1,2,3 \} \hbox{ for all }i = 0, \ldots, n-1 \}\\
C &= \{ 0 \leq c < M^n: d_i(c) \in \{ 2,3,4,5 \} \hbox{ for all }i = 0, \ldots, n-1 \}\\
D &= \{ 0 \leq d < M^n: d_i(d) \in \{ 4,5,6,8 \} \hbox{ for all }i = 0, \ldots, n-1 \}\\
G &= \{ 0 \leq a,b < M^n: (d_i(a),d_i(b)) \in \{ (4,0), (2,1), (3,1), (4,1), (0,2), (2,2), (0,3), (2,3) \}\\
& \hbox{ for all }i = 0, \ldots, n-1 \}
\end{align*}
and $M > 8$.

The authors thank Jean Bourgain for helpful discussions and for suggesting the use of more than three slices, and Christoph
Thiele for bringing the authors together.  The authors also thank Keith Rogers for pointing out some minor errors in an earlier version.  Part of this work was conducted at
ANU, UNSW, and ETH.  The 
authors are supported by NSF grants DMS-9801410 and DMS-9706764 
respectively.

\section{A combinatorial lemma}

We shall need the following combinatorial lemma.

\begin{lemma}  Let $X$ and $A_1, \ldots, A_n$ be finite sets for some
$n \geq 0$, and
for each $1 \leq i \leq n$ let $f_i: X \to A_i$ be a function.  Then
\be{iter} 
\# \{ (x_0, \ldots, x_n) \in X^{n+1}: f_i(x_{i-1}) = f_i(x_i) \hbox{ for all } 0 \leq i \leq n\}
\geq \frac{(\# X)^{n+1}}{\prod_{i=1}^n \# A_i}.
\end{equation}
\end{lemma}

The reader may verify from probabilistic methods that \eqref{iter} is sharp.

\begin{proof}  We prove by induction on $n$.  The claim is trivial for $n=0$.
Now suppose that $n \geq 1$, and the claim has been proven for $n-1$.

Define an element $a \in A_n$ to be \emph{popular} if
$$ \# \{ x \in X: f_n(x) = a \} \geq \frac{\# X}{2 \# A_n},$$
and define $X'$ to be those elements $x \in X$ such that
$f_n(x)$ is popular.  Since each unpopular element of $A_n$ contributes at most
$\# X/(2 \# A_n)$ elements to $X$, we see that
$$ \# (X \backslash X') \leq \frac{1}{2} \# X,$$
so
\be{x-card} \# X' \geq \frac{1}{2} \# X.
\end{equation}
By applying the induction hypothesis to $X'$ we have
$$ 
\# \{ (x_0, \ldots, x_{n-1}) \in (X')^n: f_i(x_{i-1}) = f_i(x_i) \hbox{ for all }0 \leq i < n\} \geq \frac{(\# X')^n}{\prod_{i=1}^{n-1} \# A_i}.$$
Since $f_n(x_{n-1})$ is popular, we thus have
$$ 
\# \{ (x_0, \ldots, x_{n-1},x_n) \in (X')^n \times X: f_i(x_{i-1}) = f_i(x_i) \hbox{ for all } 0 \leq i \leq n\} \geq \frac{(\# X')^n}{\prod_{i=1}^{n-1} \# A_i} \frac{\# X}{2 \# A_n}.$$
>From \eqref{x-card} we thus have
\be{cool} \hbox{LHS of \eqref{iter}} \geq 2^{-n} \frac{(\# X)^{n+1}}{\prod_{i=1}^n \# A_i}.
\end{equation}
To eliminate the factor of $2^{-n}$, we let $M$ be a large integer, and apply \eqref{cool} with  $X$, $A_i$ replaced by $X^M$, $A_i^M$, and $f_i$ replaced with the function $f_i^M: X^M \to A_i^M$ defined by
$$ f_i^M(x^1, \ldots, x^M) = (f_i(x^1), \ldots, f_i(x^M)),$$
to obtain
$$
(\hbox{LHS of \eqref{iter}})^M \leq 2^{-n} \left( \frac{(\# X)^{n+1}}{\prod_{i=1}^n \# A_i} \right)^M.
$$
The claim then follows by letting $M \to \infty$ (cf. \cite{ruzsa}).
\end{proof}

\section{Proof of \eqref{diff-6}}

Fix $A$, $B$, $C$, $N$.
By removing redundant elements of $G$, we may assume
\be{one}
\hbox{The map } (a,b) \mapsto a-b \hbox{ is one-to-one on } G,
\end{equation}
in which case we need to show
\be{11-6}
\# G \leq N^{11/6}.
\end{equation}

Define the set
\be{vert}
V = \{(a,b,b') \in A \times B \times B: (a,b), (a,b') \in G\}.
\end{equation}
By applying \eqref{iter} with $n=1$ and $f_1: G \to A$ being the
projection map, we see from \eqref{ab-card} that
\be{v-card}
\# V \geq \frac{(\# G)^2}{N}.
\end{equation}
Now consider the maps $f_1: V \to C \times C$, $f_2: V \to B \times B$,
$f_3: V \to C \times B$ by
\bas
f_1(a,b,b') &= (a+b,a+b')\\
f_2(a,b,b') &= (b,b')\\
f_3(a,b,b') &= (a+b,b').
\end{align*}
Let $S$ denote the set
$$ S = \{ (v_0,v_1,v_2,v_3): f_1(v_0) = f_1(v_1), f_2(v_1) = f_2(v_2), f_3(v_2) = f_3(v_3) \};$$
from \eqref{iter}, \eqref{ab-card} and eqref{c-card} we have
\be{s-card}
\# S \geq \frac{(\# V)^4}{N^6}.
\end{equation}
Write $v_i = (a_i, b_i, b'_i)$ for $i=0,1,2,3$, and consider the map $g: S \to V \times A \times B$ defined by
$$ g(v_0,v_1,v_2,v_3) = (v_0, a_2, b_3).$$
We now observe that $g$ is injective, or in other words that $(v_0,a_2,b_3)$
determines $(v_0,v_1,v_2,v_3)$.  To see this, we note from construction
\be{first}
a_0 + b_0 = a_1 + b_1, \quad a_0 + b'_0 = a_1 + b'_1, \quad b_1 = b_2, \quad b_2 = b'_2
\end{equation}
and
\be{second}
a_2 + b_2 = a_3 + b_3, \quad b'_2 = b'_3.
\end{equation}
>From \eqref{first} we have
$$ b_0 - b'_0 = b_1 - b'_1 = b_2 - b'_2$$
while from \eqref{second} we have
$$ a_3 - b'_3 = a_3 + b_3 - b_3 - b'_2 = a_2 + b_2 - b'_2 - b_3.$$
Combining these equations, we obtain
$$ a_3 - b'_3 = a_2 - b_3 + b_0 - b'_0.$$
Thus $a_3 - b'_3$ is determined by $(v_0,a_2,b_3)$.  Since $(a_3,b'_3) \in G$, we see from \eqref{one} that $(a_3,b'_3)$, and thus $v_3$, are determined by $(v_0,a_2,b_3)$.  The
injectivity of $g$ then follows by using \eqref{second} to determine $v_2$, and then \eqref{first} to determine $v_1$.

Since $g$ is injective, we have from \eqref{ab-card} that
$$ \# S \leq N^2 \# V.$$
Combining this with \eqref{s-card} we see that
$$ \# V \leq N^{8/3},$$
and \eqref{11-6} follows from \eqref{v-card}.

\section{Proof of \eqref{diff-4}}

Fix $A$, $B$, $C$, $D$, $N$.  We may assume \eqref{one} as before, so that we need to show that
\be{targ-4} \# G \leq N^{7/4}.
\end{equation}

Define $V$ by \eqref{vert} as before, and define the function $f_4: V \to D \times B$ by
$$ f_4(a,b,b') = (a+2b, b').$$
Let $T$ denote the set
$$ T = \{ (v_0,v_1) \in V^2: f_4(v_0) = f_4(v_1)\}.$$
>From \eqref{iter}, \eqref{ab-card}, \eqref{d-card} we have 
\be{t-card}
\# T \geq \frac{(\# V)^2}{N^2}.
\end{equation}
Write $v_i = (a_i,b_i,b'_i)$ for $i=0,1$, and define the map
$h: T \to C \times C \times B$ by
$$ h(v_0,v_1) = (a_0 + b_0, a_0 + b'_0, b_1).$$
We claim that $h$ is injective, so that $h(v_0,v_1)$
determines $(v_0,v_1)$.  From construction we have
\be{third}
a_0 + 2b_0 = a_1 + 2b_1, \quad b'_0 = b'_1,
\end{equation}
so that
$$ a_1 - b'_1 = a_0 + a_1 + 2b_1 - 2b_1 - (a_0 + b'_1) = 2(a_0 + b_0) - 2b_1 - (a_0 + b'_0).$$
Thus $a_1 - b'_1$ is determined by $h(v_0,v_1)$.  From \eqref{one} we thus
see that $(a_1,b'_1)$ is determined by $h(v_0,v_1)$.  The injectivity of $h$ then follows from \eqref{third}.

Since $h$ is injective, we see from \eqref{ab-card} and \eqref{c-card} that
$$ \# T \leq N^3.$$
>From \eqref{t-card} we thus have
$$ \# V \leq N^{5/2},$$
and \eqref{targ-4} follows from \eqref{v-card}.

\end{document}